\documentclass[12pt]{article}
\usepackage{amssymb}
\usepackage{graphicx}
\usepackage{latexsym}
\def\part#1{\frac{\partial\phantom{q}}{\partial#1}}

\newenvironment{rmk}{\begin{trivlist}\item[]{\bf Remark:} }
{\end{trivlist}}

\newenvironment{prf}{\begin{trivlist}\item[]{\bf Proof:} }
{\hfill $\Box$ \end{trivlist}}

\newtheorem{thm}{Theorem}

\newtheorem{prp}[thm]{Proposition}

\newcommand{\lie}[1]{\mathfrak{#1}}

\def\Hom{\mathop{\rm Hom}\nolimits}

\def\ker{\mathop{\rm ker}\nolimits}

\def\deg{\mathop{\rm deg}\nolimits}

\newcommand{\C}{\mathbf{C}}

\newcommand{\Z}{\mathbf{Z}}

\newcommand{\RP}{{\mathbf R}{\rm P}}

\begin{document}
\title{Vector bundles and the icosahedron}
 \author{Nigel Hitchin\\[5pt]}
\maketitle
\centerline{\it Dedicated to S. Ramanan on the occasion of his 70th birthday}

\section{Introduction}
 
Describing the equation of a plane curve $C$ as the determinant  of a matrix of linear forms is a classical problem. A related issue is expressing $C$ as the Pfaffian of a skew-symmetric matrix of forms. When the curve is defined by a harmonic homogeneous polynomial (a solution to  Laplace's equation), then the action of the Lie algebra of the orthogonal group $SO(3,\C)$ provides a natural such expression, and with it, a natural semi-stable rank two vector bundle $E$ on $C$. This paper concerns itself with these bundles, which have  trivial determinant and $H^0(C,E)=0$. 

For a cubic curve,  Atiyah's classification shows that  they fall into three types: a sum of line bundles $L\oplus L^*$ with $L^2$ nontrivial, a non-trivial  extension of a line bundle of order two by itself, and a trivial extension. We shall show  that  the second and third cases are described in the space of harmonic cubics by the vanishing of a certain invariant polynomial  of degree six.   

Our route to this result  uses  the Mukai-Umemura threefold and we link up with the classical  geometry of the Clebsch diagonal cubic  surface, certain distinguished rational curves on it, and the action of the symmetric group $S_5$.

The original motivation for this paper comes from the author's interest in explicit solutions to Painlev\'e equations, as in  \cite{Hit1} and \cite{Hit2}. These were connected to the study of certain threefolds with an open orbit of $SO(3,\C)$ and finite stabilizer.  The Mukai-Umemura manifold has such an action where the stabilizer is the icosahedral group. Since the first paper, where the stabilizer is the dihedral group,  was written for the 60th birthdays of Narasimhan and Seshadri, it seems appropriate, sixteen  years on, to discuss the icosahedron for Ramanan. 

\section{Representations of $SO(3,\C)$}
Let $V$ be a finite-dimensional  irreducible representation space of the  group $SO(3,\C)$. These occur in each odd dimension $(2d+1)$ and have an invariant inner product $(u,v)$. Restricted to the circle subgroup $SO(2,\C)\subset SO(3,\C)$, the weights are $-d\le n\le d$;  in particular the weight zero occurs with multiplicity one so there is a unique invariant element. The usual realization of $V$ is as $S^{2d}$ -- the space of homogeneous polynomials in the complex variables $z_1,z_2$ of degree $2d$ under the action of $SL(2,\C)/\pm 1\cong SO(3,\C)$.

The Lie algebra $\lie{so}(3)$ acts  on $V$ as skew-adjoint transformations and so each $x\in \lie{so}(3)$ defines a skew form 
$$\omega_x(u,v)=(x\cdot u,v).$$
Now fix $v\in V$ and restrict $\omega_x$ to the $2d$-dimensional orthogonal complement $W$ of $v$. The  Pfaffian $\omega_x\wedge\omega_x\dots\wedge \omega_x \in \Lambda^{2d}W^*\cong \C$ defines a homogeneous polynomial $f(x)$ of degree $d$. 

 The fact that we have a natural map from vectors in $V$ to polynomials is no surprise if we recall that another  realization of the $(2d+1)$-dimensional irreducible representation of $SO(3,\C)$ is as spherical harmonics -- homogeneous polynomials  $f(x_1,x_2,x_3)$ of degree $d$ which satisfy Laplace's equation. However, $v$ gives somewhat more than just a polynomial. 
 
 The skew form $\omega_x$ defines a map $\omega_x: W\rightarrow W^*$, which is linear in $x$. Hence if ${\mathbf P}^2$ is the projective space of the vector space $\lie{so}(3)$, we have a sequence of sheaves

\begin{equation}
0\rightarrow {\mathcal O}_{{\mathbf P}^2}(W(-2))\stackrel{\omega} \rightarrow {\mathcal O}_{{\mathbf P}^2}(W^*(-1))\rightarrow E\rightarrow 0
\label{exact1}
\end{equation}

where $E$ is a rank $2$ vector bundle supported on the curve $C\subset {\mathbf P}^2$ defined by $f(x)=0$. The vector bundle satisfies $\Lambda^2E\cong K_C$, the canonical bundle, and from the exact cohomology sequence, $H^0(C,E)=0$. 

These facts about curves defined by Pfaffians (and much more besides) can be found in \cite{Bea}. In particular, the bundle is always semi-stable. To see this, note that  if $L\subset E$ is a subbundle,  the inclusion defines a section $i$ of $L^*E$ and if $\deg L\ge g$ then by Riemann-Roch there is a non-zero section $s$ of $L$. But then $si$ is a nonzero section of $E$ which is a contradiction. Hence we must have
$$\deg L\le g-1=\frac{1}{2}\deg K_C=\frac{1}{2}\deg E.$$

This natural process thus generates a curve and a rank $2$ semi-stable bundle, and the question we ask is ``What is this bundle?"
\vskip .5cm
When $d=3$, the curve $C$ is a plane cubic and then we know from Atiyah's classification of bundles on an elliptic curve \cite{At} that, when $C$ is smooth, there are only three possibilities.  Since $K_C\cong {\mathcal O}$,  $\Lambda^2E$ is trivial and $H^0(C,E)=0$ and we have the cases: 
\begin{itemize}
\item
$E=L\oplus L^*$ where $\deg L=0$ and $L\ne {\mathcal O}$ (the generic case) 
\item
$E$ is a non-trivial extension $0\rightarrow L\rightarrow E\rightarrow L\rightarrow 0$
 where $L^2\cong {\mathcal O}$ 
 and $L\ne {\mathcal O}$ 
\item
$E=L\oplus L$ where  $L^2\cong {\mathcal O}$ 
 and $L\ne {\mathcal O}.$
\end{itemize}
These three types define a stratification of the six-dimensional projective space ${\mathbf P}(V)$ where $V$ is the seven-dimensional space of harmonic cubic polynomials $f$ --  or at least that part of it for which  the curve $f(x)=0$ is nonsingular.   To investigate this further we have to see a link with another piece of geometry.

\section{Isotropic spaces}\label{iso}

Let $d=3$, and as above let $f\in V$ define a smooth cubic curve $C\subset {\mathbf P}^2$ and  $W\subset V$ the six-dimensional orthogonal complement of $f$ in $W$.  Since $\Lambda^2E\cong {\mathcal O}$, we have $E^*\cong E$ and then the kernel of $\omega:W\rightarrow W^*(1)$ on $C$ can be identified with $E(-1)$. 

Suppose $U\subset W$ is a three-dimensional subspace, {\it isotropic} with respect to the skew forms $\omega_x$ for all $x$. Then, for each $x$, $\omega_x:W\rightarrow W^*$ maps $U$ to its annihilator $U^o$ and so defines a homomorphism of sheaves on ${\mathbf P}^2$, $\omega\vert_U:U\rightarrow U^o(1)$. Moreover the determinant of this map is the Pfaffian of $\omega_x$. The kernel of $\omega\vert_U$ on $C$ is then a  line bundle $L(-1)\subset U$. Since this is the restriction of $\omega$ to $U$, $L$ is a subbundle of $E$.  From the treatment of determinantal loci in \cite{Bea} we have $\deg L=0$ and $H^0(C,L)=0$ so $L\ne {\mathcal O}$. 

Conversely, on $C$,   
suppose we have a line bundle $L\subset E$ with $\deg L=0$ and $L\ne {\mathcal O}$.  From the exact cohomology sequence of (\ref{exact1}) we have $W^*\cong H^0(C,E(1))$. If  $\alpha\in W^*$ annihilates $L(-1)$ then it defines a section of  $(E/L)(1)=L^*(1)$ which is of degree $3$ and so has a three-dimensional space of sections on $C$. It follows that $L(-1)\subset W$ sweeps out a three-dimensional subspace $U$.

 Now $\omega\vert_U$ is a section of  $\Lambda^2U^*(1)$ on $C$ and we have 
$\Lambda^2U^*\cong U\otimes\Lambda^3U^*$ so this defines a section $w$ of $U(1)$ which tautologically is annihilated by $\omega$.
Thus if $w\ne 0$, then $\ker \omega\vert_U \cong {\mathcal O}(-1)$.  But $\ker \omega\vert_U\cong L(-1)$ and $L\ne {\mathcal O}$, so we deduce that $w=0$ and $U$ is isotropic for each $\omega_x$.
 
 If we look again at the trichotomy of Atiyah's classification,  we see that the three cases are equivalent to
 \begin{itemize}
\item
$E$ has two rank one subbundles 
\item
$E$ has one rank one subbundle 
\item
$E$ has infinitely many rank one subbundles 
\end{itemize}
and by what we have just seen, this condition translates into an equivalent statement about the three-dimensional  isotropic subspaces of $W$:
 \begin{itemize}
\item
$W$ has two  isotropic subspaces 
\item
$W$ has one  isotropic subspace  
\item
$W$ has infinitely many isotropic subspaces.
\end{itemize}

 This provides an extension of the criterion to {\it any} $f\in V$, and not just those which define smooth cubic curves. (In fact, Atiyah's classification has been extended in \cite{Burb} to a class of singular elliptic curves; moreover the wild case of cuspidal cubics does not occur when $f$ is harmonic).
 
\section{The Mukai-Umemura threefold}

The study of three-dimensional isotropic subspaces of $V$ is best approached via a special Fano threefold introduced by Mukai and Umemura \cite{Muk}. We consider the Grassmannian $G(3,V)$ of three-dimensional subspaces of the seven-dimensional representation space $V$ and its universal rank $3$ bundle ${\mathcal E}$.  For $x\in {\lie {so}}(3)$ we have the skew form $\omega_x$ on $V$, and so a section of the rank $9$ vector bundle over $G(3,V)$:
$$\Lambda^2{\mathcal E}^*\otimes  {\lie {so}}(3).$$
Now  $\dim G(3,V)=3\times (7-3)=12$ and then the zero set of the section is a smooth $12-9=3$-dimensional manifold known as the {\it Mukai-Umemura threefold} $Z$. By construction it parametrizes subspaces $U$ isotropic for all $\omega_x$ and thereby has a natural action of $SO(3,\C)$. 

The threefold $Z$ has the same additive integral cohomology as ${\mathbf P}^3$ and $H^2(Z,\Z)$ is generated by $x=c_1$. However $x^2\in H^4(Z,\Z)$ is $22$ times a generator $y$. Since $c_1>0$ the Todd genus is $1$ so $c_1c_2=24$ and we have Chern classes
\begin{equation}
c_1=x \qquad c_2=24y \qquad c_3=4xy.
\label{chern}
\end{equation}

Given $f\in V$, the inner product $(f,-)$ defines an element of $V^*$ and, restricting to the universal bundle ${\mathcal E}\subset V$, we get a section of ${\mathcal E}^*$. Over $Z$, this is a rank three bundle and the section vanishes at the points which correspond to isotropic three-dimensional subspaces orthogonal to $f$. 

\begin{prp} On the Mukai-Umemura threefold, $c_3({\mathcal E}^*)=2$.
\end{prp}
\begin{prf} The tangent bundle of the Grassmannian is $\Hom ({\mathcal E}, V/{\mathcal E})$ and $Z$ is the non-degenerate zero set of a section of $\Lambda^2{\mathcal E}^*\otimes  \C^3.$ Thus, as $C^{\infty}$ bundles
$$TZ\oplus (\Lambda^2{\mathcal E}^*\otimes  \C^3)\cong \Hom ({\mathcal E}, V/{\mathcal E}).$$
Applying the Chern character we find
$$c_1({\mathcal E}^*)=c_1\qquad c_2({\mathcal E}^*)=c_1^2-\frac{1}{2}c_2\qquad c_3({\mathcal E}^*)=\frac{1}{10}(c_3+4c_1^3-3c_1c_2).$$
From (\ref{chern}) we obtain $10c_3({\mathcal E}^*)= 4+88-72=20$ and hence the result.
\end{prf}

A section with nondegenerate zero set will thus vanish at two points.  We can now see the trichotomy in terms of the section $s$ of ${\mathcal E}^*$:
\begin{itemize}
\item
$s$ vanishes at two points
\item
$s$ vanishes at one point
\item
$s$ vanishes on a subvariety of positive dimension. 
\end{itemize}

\section{The icosahedron}
The ten-dimensional space of all homogenous cubics in $x$ has an $SO(3,\C)$-invariant product which can be normalized so that $(f(x),(x,a)^3)=f(a)$. For fixed $a\in \C^3$, the polynomial $(x,a)^3$ is not harmonic but  its orthogonal projection onto $V$ is
$$f_a(x)=(x,a)^3-\frac{3}{5}(a,a)(x,a)(x,x).$$ 
This is invariant by the group of rotations fixing $a$. For any cubic $f$ we still have the inner product property $(f,f_a)=f(a)$.

The action of $u\in {\lie {so}}(3)$ on $f_a$ is
$$(u\cdot f_a)(x)=3(x,[u,a])[(x,a)^2-\frac{1}{5}(a,a)(x,x)].$$

Take a regular icosahedron with vertices at $\pm a_1,\pm a_2\dots,\pm a_6$. The angle $\theta$ between any two axes joining opposite vertices satisfies $\cos \theta =\pm 1/\sqrt{5}$. Now consider 
$$(u\cdot f_{a_i},f_{a_j})=3(a_j,[u,a_i])[(a_j,a_i)^2-\frac{1}{5}(a_i,a_i)(a_j,a_j)].$$
If $i=j$ then $(a_j,[u,a_i])=0$ and when $i\ne j$ 
$$(a_j,a_i)^2-\frac{1}{5}(a_i,a_i)(a_j,a_j)=(\cos^2 \theta-\frac{1}{5})(a_i,a_i)(a_j,a_j)=0$$
so $f_{a_1},\dots, f_{a_6}$ span an isotropic subspace $U$, invariant under the icosahedral group $G$. Its orthogonal complement in $V$ consists of those $f$ such that 
$$0=(f,f_{a_i})=f(a_i)=0.$$

There are five objects which are permuted by $G$ which realize the well-known isomorphism $G\cong A_5$. These are the five sets of three orthogonal planes in which all $12$ vertices lie. Let $e_1,e_2,e_3$ be the three unit normals to such a set  then the cubic 
$(e_1,x)(e_2,x)(e_3,x)$ vanishes at $a_i$ as do all its transforms by $G$. These  span the four-dimensional permutation representation  ${\bf 4}$ of $A_5$. Moreover, since   $(e_1,x)(e_2,x)(e_3,x)$ satisfies Laplace's equation, this is the orthogonal complement of $U$ in $V$, hence $U$ has dimension $3$.
 Every icosahedron therefore defines a point in the Mukai-Umemura manifold and the three-dimensional orbit $SO(3,\C)/G$ is a dense open set in $Z$. 

Generically, a harmonic cubic $f$ defines a section $s$ of ${\mathcal E}^*$ which vanishes at two points in this open set, which means that the corresponding plane cubic $C$ contains the twelve points in ${\mathbf P}^2$ defined by the axes of two icosahedra. Put another way, if a generic cubic curve contains  one such ``icosahedral set"  $\{[a_1],[a_2],\dots,[a_6]\}\subset {\mathbf P}^2$ then it contains another.  

\section{Degenerate icosahedra}
\subsection{The two types}

We need to understand also the divisor of ``degenerate icosahedral sets" which forms the complement of the open orbit of $SO(3,\C)$  in $Z$. These correspond to isotropic subspaces of two types, constituting two orbits, of dimension $2$ and $1$ respectively. 

A representative of the first type is the space with basis 
$$(b,x)^3-\frac{3}{5}(b,x)(b,b)(x,x),\quad (a,x)^2(b,x),\quad (a,x)^3$$
where $(a,a)=0, (a,b)=0$ and $(b,b)\ne 0$. This subspace is invariant by the rotations about the axis $b$ and is spanned by the weight spaces $\{0,2,3\}$ for that action. Geometrically, the isotropic subspace is defined by $[b]\in {\mathbf P}^2$ (a point not on  the null conic $Q$ defined by $(x,x)=0$) together with  the point of contact $[a]$ of a tangent to $Q$ through $[b]$. The choice of tangent is a double covering of ${\mathbf P}^2\setminus Q$, which  is abstractly an affine quadric. 

The icosahedral set here describes five of the vectors   $\{a_1,a_2,\dots,a_6\}$ coalescing into a single  vector $a$. Under this degeneration the  relation $(a_j,a_i)^2-\frac{1}{5}(a_i,a_i)(a_j,a_j)=0$ implies that $a$ is null and the sixth vector $b$ satisfies $(b,a)=0$.

The second type is a subspace with basis 
$$(a,x)((b,x)^2-\frac{1}{5}(b,b)(x,x)),\quad  (a,x)^2(b,x),\quad (a,x)^3$$ 
where $(a,a)=0$ and $(a,b)=0$. This is geometrically defined by the point $[a]\in Q$ and is invariant by the Borel subgroup which fixes $[a]$. Decomposing with respect to a semisimple element this is the span of weight spaces   $\{1,2,3\}$. Here all six  vectors   $\{a_1,a_2,\dots,a_6\}$ coalesce into a single null vector $a$.

The union of these two orbits forms an anticanonical divisor in $Z$, for if $X_1,X_2,X_3$ are the vector fields on $Z$ generated by a basis of ${\lie {so}}(3)$, then $X_1\wedge X_2\wedge X_3$ vanishes on the lower-dimensional  orbits and this is a section of $K_Z^*=\Lambda^3TZ$.  Since $c_1(Z)$ is a generator of $H^2(Z,\Z)$ it must vanish with multiplicity $1$. Note that this anticanonical divisor $D$ cannot be smooth (for then it would be a K3 surface or a torus). It is instead a singular image of a map $\alpha:{\mathbf P}^1\times {\mathbf P}^1\rightarrow Z$ (in fact, a transverse to  the diagonal in ${\mathbf P}^1\times {\mathbf P}^1$ has a cusp singularity $y^2=x^3$ (see \cite{SKD}).

To define $\alpha$ recall that we have identified the type 1 orbit as an affine quadric, the complement of the diagonal in ${\mathbf P}^1\times {\mathbf P}^1$. This is  an ordered pair of points $[a],[a']$ on the conic $Q\cong {\mathbf P}^1$, where the two tangents meet at $[b]\in {\mathbf P}^2$.  To extend the map $\alpha$ to the diagonal 
put $c=a+tb$, with $(c,a)=0$ and let  $t\rightarrow 0$. Then 
$$(c,x)^3-\frac{3}{5}(c,x)(c,c)(x,x)=(a,x)^3+3t(a,x)^2(b,x)+3t^2[(a,x)((b,x)^2-\frac{1}{5}(b,b)(x,x))]. .$$
and the three leading coefficients span  the type 2 isotropic subspace.

 Consider a degenerate icosahedral set of type 1, spanned by the three cubic polynomials  $(b,x)^3-{3}(b,x)(b,b)(x,x)/5, (a,x)^2(b,x)$ and $(a,x)^3$. Then  $f$ is orthogonal to this if the curve $C$ given by $f(x)=0$ intersects the conic $Q$ tangentially at $[a]$ and passes through $[b]$. For type 2, the intersection multiplicity of $C$ with $Q$ at $[a]$  must be at least $3$.
 
 \subsection{The universal bundle}
 
 We should also consider the universal bundle ${\mathcal E}$ on the divisor $D$, or rather its pullback $\alpha^*{\mathcal E}$ on ${\mathbf P}^1\times {\mathbf P}^1$. Note that both types of degenerate isotropic subspace contain $(a,x)^3$ and $(a,x)^2(b,x)$ where $(b,a)=0$, or equivalently the two-dimensional subspace given by  $(a,x)^2(c,x)$ for all $c$ orthogonal to $a$. 
 
 We use the two-fold covering map $\pi: {\mathbf P}^1\times {\mathbf P}^1\rightarrow {\mathbf P}^2$, the quotient by the involution interchanging the factors. Each factor is isomorphic to the diagonal which maps to the conic $Q$. As usual, geometrically $\pi([a],[a'])$ is the point $[b]$ of intersection of the tangents at $[a],[a']\in Q$.

 The two-dimensional vector space $a^{\perp}\subset \C^3$ for $[a]\in Q$  defines a vector bundle $A$ over $Q$ for which the projective bundle ${\mathbf P}(A)$ is trivial. Indeed,  it is the bundle of tangent lines to the conic and so under the map $\pi$ can be identified with the family ${\mathbf P}^1\times \{x\}$, $x\in {\mathbf P^1}$, which is a trivial bundle over the diagonal. On the other hand since $(a,a)=0$ 
 $$a^{\perp}\cong (\C^3/{\mathcal O}_{{\mathbf P}^2}(-1))^*$$
 and so $\Lambda^2 A \cong {\mathcal O}_{{\mathbf P}^2}(-1)={\mathcal O}(-2)$, identifying $Q$ with ${\mathbf P}^1$. Hence, since  ${\mathbf P}(A)$ is trivial,
 $$A\cong {\mathcal O}(-1)\otimes \C^2.$$
 Multiplying by the factor $(a,x)^2$ which is quadratic in $a$, it follows that the subbundle in $\alpha^*{\mathcal E}$ of cubics of the form $(a,x)^2(c,x)$ with $(c,a)=0$   is isomorphic to ${\mathcal O}(-5,0)\otimes \C^2.$ 
 
The cubic $(b,x)^3-3(b,x)(b,b)(x,x)/5$, together with the subspace $A$, spans  a type 1 degenerate subspace  for $(b,b)\ne 0$ and when $(b,b)=0$ it lies in $A$. This term is homogeneous of degree $3$ in $b$, and so defines a homomorphism from $\pi^*{\mathcal O}_{{\mathbf P}_2}(-3)={\mathcal O}(-3,-3)$ to $\alpha^*{\mathcal E}$. It projects to a section of the line bundle  $\alpha^*{\mathcal E}/A$ which vanishes  on  the diagonal $\Delta\subset {\mathbf P}^1\times {\mathbf P}^1$. Hence $\alpha^*{\mathcal E}/A$ is of the form ${\mathcal O}(k,k)$ and 
$$\Lambda^3\alpha^*{\mathcal E}\cong  {\mathcal O}(k-10,k).$$

But for $Z$ we have  the Chern number $c_1^3=22$ and since $D$ is an anticanonical divisor, $c_1^2[D]=22$. However $c_1({\mathcal E})=-c_1(Z)$ and so 
$$22=2k(k-10)$$
and $k=-1$. Thus $\alpha^*{\mathcal E}$ is an extension
$$0\rightarrow {\mathcal O}(-5,0)\otimes \C^2\rightarrow \alpha^*{\mathcal E}\rightarrow {\mathcal O}(-1,-1)\rightarrow 0.$$
 This extension is classified by an element of $H^1( {\mathbf P}^1\times {\mathbf P}^1,{\mathcal O}(-4,1))\otimes \C^2$. It is also by definition $SO(3,\C)$-invariant by the diagonal action on ${\mathbf P}^1\times {\mathbf P}^1$. But as representation spaces 
 $$H^1( {\mathbf P}^1\times {\mathbf P}^1,{\mathcal O}(-4,1))\otimes \C^2\cong H^1({\mathbf P}^1,{\mathcal O}(-4))\otimes  H^0({\mathbf P}^1,{\mathcal O}(1))\otimes \C^2\cong S^2\otimes S\otimes S$$
  and $\alpha^*{\mathcal E}$ is defined by a non-zero vector in the unique invariant one-dimensional subspace. 
  
  The sections of $\alpha^*{\mathcal E}^*$ fit into the exact sequence
  $$0\rightarrow H^0( {\mathbf P}^1\times {\mathbf P}^1,{\mathcal O}(1,1))\rightarrow  H^0( {\mathbf P}^1\times {\mathbf P}^1,\alpha^*{\mathcal E^*})\rightarrow H^0({\mathbf P}^1,{\mathcal O}(5))\otimes \C^2\rightarrow 0.$$
  As a representation space we then have 
  $$H^0( {\mathbf P}^1\times {\mathbf P}^1,\alpha^*{\mathcal E^*})\cong S\otimes S+S\otimes S^5$$
 which contains with multiplicity one the seven-dimensional representation $S^6$ as a subspace of $S\otimes S^5$. These sections are the restriction to $D$ of the sections $(f,-)$ on $Z$ that we have considered earlier. There is one important consequence of this:
 \begin{prp} \label{ten} A harmonic cubic $f$ is orthogonal to a degenerate isotropic subspace if and only if $\Delta(f)=0$ for a certain $SO(3,\C)$-invariant polynomial $\Delta$ of degree $10$.
 \end{prp}
 \begin{prf} We want to know when a section $s$ of ${\mathcal E^*}$ defined by $(f,-)$ vanishes on $D$. If it does vanish somewhere then the map to $H^0({\mathbf P}^1,{\mathcal O}(5))\otimes \C^2$ gives two sections of ${\mathcal O}(5)$ with a common zero. These arise from $S^6\subset S^5\otimes S$ so the condition is that we have a homogeneous polynomial $p(z_1,z_2)$ of degree $6$ such that the two partial derivatives $\partial p/\partial z_1,\partial p/\partial z_2$, homogeneous of degree $5$, have a common zero. The vanishing of the resultant is the condition. This is a degree $10$ polynomial $\Delta$ in the coefficients of $f$, the {\it discriminant}. Its vanishing  implies that there is a point $[a]=[z_1,z_2]\in {\mathbf P}^1$ where the section  of ${\mathcal O}(6)$ has a double zero. Then the section of ${\mathcal O}(5,0)\otimes \C^2$ vanishes on $\{[a]\}\times  {\mathbf P}^1\subset {\mathbf P}^1\times {\mathbf P}^1$.
 
 The dual of the inclusion ${\mathcal O}(-3,-3)\subset \alpha^* {\mathcal E}$ gives a homomorphism  
 $$\alpha^*{\mathcal E}^*\rightarrow {\mathcal O}(3,3)$$
 which determines the third component of the section $s$. It maps $H^0( {\mathbf P}^1\times {\mathbf P}^1,{\mathcal E^*})$ to $$H^0( {\mathbf P}^1\times {\mathbf P}^1,{\mathcal O}(3,3))\cong S^3\otimes S^3$$ and our seven-dimensional representation space $S^6$ maps into the symmetric elements. Such a section is then the pull-back of a section of  ${\mathcal O}_{{\mathbf P}^2}(3)$ - the cubic $f$. In particular, its divisor  passes through $([a],[a])$ where the other two components of $s$ vanish. 
 
 We see then that the section $s$ vanishes at a point of $D$ if the discriminant vanishes. Geometrically this implies that the cubic curve $C$ meets the null conic tangentially at $[a]$.

 \end{prf}

\section{The Clebsch cubic surface}\label{clebsch}

The standard permutation representation of $A_5$ is the  four-dimensional subspace ${\bf 4}\subset \C^5$ defined by $y_1+y_2+\dots +y_5=0$. We shall investigate next the  geometry of the corresponding three-dimensional projective space ${\mathbf P}^3$. 

The ring of invariants for $A_5$ on ${\bf 4}$ is generated by the elementary symmetric functions $\sigma_2,\sigma_3,\sigma_4,\sigma_5$ in $(y_1,y_2,\dots,y_5)$ and the degree $10$ invariant
$$\prod_{i<j}(y_i-y_j).$$
The invariant $y_1^2+y_2^2+\dots +y_5^2$ defines a nonsingular quadratic form on ${\bf 4}$ and then  orthogonality of subspaces defines polarity in ${\mathbf P}^3$: each point has a polar plane, and each line  a polar line. There is a {\it unique} invariant cubic, which we can write as either $\sigma_3=0$ or more usually as 
$$y_1^3+y_2^3+\dots +y_5^3=0.$$
This equation defines the {\it Clebsch cubic surface} $S$. It is nonsingular and is invariant under the action of the full symmetric group $S_5$. 

Consider now  six points $[a_1],[a_2],\dots,[a_6]\in {\mathbf P}^2$ forming an icosahedral set. No three are collinear and no conic passes through them all, for then it would be invariant by the icosahedral group and so given  by the null conic $(x,x)=0$; but the vertices $a_i$ are not null.  It follows that blowing up the six points gives a non-singular cubic surface in ${\mathbf P}^3$. More precisely, the embedding is given by the plane cubic curves that pass through the six points, which is the representation ${\bf 4}$ of $A_5$. Hence, by uniqueness, the invariant cubic surface  must be the Clebsch surface $S$. 

The blown up points form six disjoint lines in $S$, half of a double-six configuration. The other six are given by the proper transforms of conics passing through five of the six points. These we have encountered already:

$$(a_i,x)^2-\frac{1}{5}(a_i,a_i)(x,x)=0.$$

Blowing these down gives another map from $S$ to a projective plane $\tilde{\mathbf P}^2$. Moreover since the icosahedral group permutes these lines, it is the plane of a three-dimensional representation. These two planes are the projective spaces of the two inequivalent three-dimensional representations ${\bf 3}$ and $\tilde{\bf 3}$ of $A_5$. (These can be viewed as the self-dual and anti-self-dual two-forms on ${\bf 4}$).

Since the null conic $Q\subset {\mathbf P}^2$ does not meet the points $[a_i]$ to be blown up, it lifts to a rational curve $R\subset S$ with self-intersection number $4$, and hence of degree $6$ in ${\mathbf P}^3$. The null conic in the second projective plane similarly lifts to a degree $6$ curve $\tilde R\subset S$. Each of these is individually invariant by the group $A_5$ and interchanged by  $S_5$. 

The rational curves $R$ and $\tilde R$ have a rather special relationship: the polar plane of  a point  $p\in R$ is a tritangent plane to $\tilde R$ and vice-versa. Furthermore the polar line of the tangent line through $p$ is a trisecant of $\tilde R$. This classical result can be found in \cite{Mel}. It establishes a $3:3$ correspondence between the two curves.

\section{The trichotomy for ${\mathbf P}^3$}
So far we have introduced three pieces of geometry -- the vector bundle on a harmonic cubic curve, the Mukai-Umemura threefold, and the projective space ${\mathbf P}^3$ of the permutation representation of $A_5$. We seek the distinguished $SO(3,\C)$-invariant subvarieties in the six-dimensional projective space ${\mathbf P}(V)$ which describe the three different types of bundles. To do this we first restrict to  ${\mathbf P}^3\subset {\mathbf P}(V)$, where we shall see the classical geometry described in the previous section playing a role.  The first step is to establish a rational correspondence between the Mukai-Umemura threefold and ${\mathbf P}^3$.

\subsection{From  Mukai-Umemura to ${\mathbf P}^3$}\label{fromMuk}

The Mukai-Umemura threefold parametrizes three-dimensional isotropic subspaces of $V$. Fix a nondegenerate one, $U_0$, defined by the icosahedral set $[a_1],[a_2],\dots,[a_6]$ as a base point. For any other $U$ we can consider the subspace $U+ U_0\subset V$. Generically this has dimension $6$ and so has a one-dimensional  orthogonal complement  spanned by a harmonic cubic $f$:  $f$ then has  two orthogonal isotropic subspaces, $U$ and $U_0$. This establishes a rational map from $Z$ to ${\mathbf P}^3={\mathbf P}(U_0^{\perp})$, but we must examine the indeterminacy.

\begin{prp}  \label{pencil} If $\dim(U+U_0)<6$ then  it is defined by an icosahedral set with a point in common with $[a_1],[a_2],\dots,[a_6]$.
\end{prp}

\begin{prf}
Suppose  that $\dim(U+U_0)<6$, then $\dim (U+U_0)^{\perp}\ge 2$ and there is a pencil of cubics passing through the two icosahedral sets. Suppose that  the  icosahedral set defined by $U$  is nondegenerate:  $[b_1],[b_2],\dots,[b_6]$. If two axes of an icosahedron coincide then so do all of the axes, so if the icosahedra are distinct there must be at least $11$ points for the pencil to pass through, 
but by  B\'ezout's theorem, unless there is a common component in the pencil, the maximum number of intersections is nine : we deduce that there must be a common component. 

If it is  a line, then since no three of the points $[a_i]$ or $[b_i]$ are collinear there is a maximum of four -- two from each set -- on the line. But then we have seven remaining for the pencil of conics. By B\'ezout again we reach a contradiction. If the common component is a conic, then it contains at most  five points from each set, and the remaining point must lie in a pencil of lines and so for some $i,j$, $[a_i]=[b_j]$.  The equation is therefore of the form 
\begin{equation}
(c,x)((a,x)^2-\frac{1}{5}(a,a)(x,x))
\label{eq2}
\end{equation}
where $[a]=[a_i]=[b_j]$ and $(c,a)=0$. 

If $f$ is orthogonal to a degenerate isotropic subspace then the cubic is tangential to the conic $Q$ at a point $[a]$ and meets the tangent  line at a point $[b]$. The pencil thus consists of the pencil of planes in ${\mathbf P}^3$ containing the tangent line to the rational curve $R\cong Q$ at the corresponding point $p$. Unless the tangent line is contained in the Clebsch cubic $S$, the generic curve in the pencil is smooth at $p$, and we can apply the Cayley-Bacharach theorem (in the degenerate case as in  \cite{GH} p.672) to deduce that  any cubic which passes  through the six  points $[b],[a_2],[a_3],\dots,[a_6]$  and is also tangent to $Q$ at $[a]$ must pass through $[a_1]$. But take the cubic which is the conic through $[a_2],[a_3],\dots,[a_6]$ together with the tangent at $[a]$:
$$(a,x)((a_1,x)^2-\frac{1}{5}(a_1,a_1)(x,x))=0.$$
If $a_1$ lies on this then $(a,a_1)=0$. Repeating for the other $[a_i]$ we deduce that they are all collinear, which is a contradiction.  

The tangents to $R$ which are contained in $S$ are the lines $\tilde E_i$, and the pencil is the same as (\ref{eq2}). 
\end{prf}

For each vertex $a_i$ of an  icosahedron,  Proposition \ref{pencil} has highlighted the role of  the icosahedra with this as a common vertex. They  form an orbit of the group of rotations fixing $a_i$ --  a  $\C^*$ orbit in $Z$. Its closure is a rational curve with  two extra points, fixed by the action, which lie in the divisor $D$ (they are  two type 1 degenerate icosahedra, determined by the two tangent lines to $Q$  from  $[b]=[a_i]$).  These six rational curves $C_i$ all pass through the basepoint $U_0$ in $Z$. From  Proposition \ref{pencil} there is a well-defined map  from  the complement of  these six curves to ${\mathbf P}^3$:
$$\beta(U)=(U+U^0)^{\perp}.$$
We shall extend this to $Z$ by blowing up certain subvarieties.

\noindent (i) First  we need to deal with the basepoint $U=U_0$. Blow up  the basepoint in $Z$, replacing it by the projectivized tangent space. Since the stabilizer of  
$U_0$ is finite, the tangent space is naturally  ${\lie{so}(3)}$. For each $b\in {\lie{so}(3)}\cong \C^3$ we orthogonally project $f_b$ onto $U_0^{\perp}$ to get $p_b$.  If this projection is zero  then $(f,f_b)=f(b)=0$ for all $f\in U_0$, so that any cubic which vanishes at the $a_{i}$ also vanishes at $b$. However, these cubics define a projective embedding of the blown-up plane, and so must separate points. It follows that  the projection is zero  if and only if $[b]=[a_{i}]$ for some $1\le i\le 6$.  We can therefore extend $\beta$ to ${\mathbf P}^2\setminus \{[a_1],[a_2],\dots,[a_6]\}$ by using the orthogonal projection $p_b$. 

\noindent (ii) The six proper transforms of the curves $C_i$ meet the exceptional divisor in the points $[a_i]$ and we need to blow up these curves  to extend the map. Consider  $p_{b(t)}$ for  $b(t)=a_i+tv$ as $t\rightarrow 0$, where $(v,a_i)=0$.  We have 
$$f_{b(t)}=f_{a_i}+3t(v,x)\left((a_i,x)^2-\frac{1}{5}(a_i,a_i)(x,x)\right)+\dots$$
 Now $f_{a_i}$ lies in $U_0$, so 
$$p_{b(t)}=3t(v,x)\left((a_i,x)^2-\frac{1}{5}(a_i,a_i)(x,x)\right)+\dots$$
and  the coefficient of $t$ extends $\beta$ to the blow-up of  the projectivized tangent space at  $[a_i]$. Its image in ${\mathbf P}^3$ is the Clebsch cubic. The exceptional curve  $E_i$ obtained  by blowing up  $[a_i]$  maps to the  pencil in Proposition \ref{pencil}. 

Denote by $\hat Z$ the blown-up Mukai-Umemura threefold.

\begin{rmk} Blowing up $C_i$ picks out a cubic of the form 
\begin{equation}
(c,x)((a_i,x)^2-\frac{1}{5}(a_i,a_i)(x,x))
\label{spec}
\end{equation}
with $(c,a_i)=0$ which extends the map $\beta$. The normal bundle of $C_i$ is trivial so the blow-up is $C_i\times {\mathbf P}^1\cong {\mathbf P}^1\times   {\mathbf P}^1$ and $\beta$ collapses the first factor. This corresponds to the fact that the cubic (\ref{spec}) contains the one-parameter family of icosahedral sets which include $[a_i]$. In general the positive-dimensional fibres of the map $\beta$ give us cubics containing  infinitely many icosahedral sets.
\end{rmk}
\subsection{The trisecant surface}

The image  $\beta(D)\subset \hat Z$ in ${\mathbf P}^3$ consists of those polynomials $f\in U_0^{\perp}$ which are also  orthogonal to a degenerate isotropic subspace.  From Proposition \ref{ten} this is given by the vanishing of an $A_5$-invariant  polynomial of degree $10$, the restriction of the $SO(3,\C)$-invariant polynomial $\Delta$.  

On the other hand, we know that if $f$   is orthogonal to a degenerate isotropic subspace, the cubic curve $C$ is defined by a plane section of $S$ which is tangent to the rational curve $R$. This means that $[f]\in {\mathbf P}^3$ lies on the polar line of a tangent line to $R$.  From the special properties of the curves $R$ and $\tilde R$ observed in Section \ref{clebsch}
  this is a {\it trisecant} of  $\tilde R$. So the degree $10$ polynomial vanishes on the trisecant surface $T$ of $\tilde R$. 
  
  Now a generic line meets  $d$ trisecants of $\tilde R$ if its polar line $\ell$ meets $d$ tangents to $R$. But this surface is the tangent developable surface  and we can calculate the degree  easily. If $X$ is a curve of degree $n$ and $\C^2\subset H^0({\mathbf P}^3,{\mathcal O}(1))$ is the pencil of planes through a line $\ell$ then under the one-jet evaluation map
  $$ev: H^0({\mathbf P}^3,{\mathcal O}(1))\rightarrow J^1(X,{\mathcal O}_{{\mathbf P}^3}(1))$$
  we have a section of 
  $$\Lambda^2(J^1(X,{\mathcal O}_{{\mathbf P}^3}(1)))\cong K_X(2)$$
  and the degree of this is $2g-2+2n$. It vanishes at the points of intersection of $\ell$ and the surface. In our case 
   $R$ is of degree $6$ and rational so the developable has degree $0-2+12=10$. 
   
   Thus the trisecant surface, which is irreducible, is of degree $10$ and given by the vanishing of the invariant polynomial $\Delta$ of  Proposition \ref{ten}. It is the locus of  cubics which pass through a nondegenerate icosahedral set $[a_1],[a_2],\dots,[a_6]$ {\it and} a degenerate one.
   
   \subsection{The three cases for ${\mathbf P}^3$}
   
   Consider now the trichotomy in the light of the construction of the map $\beta:\hat Z\rightarrow {\mathbf P}^3$. The trisecant surface $T=\beta(D)$ consists of cubics which contain a nondegenerate icosahedral set  and a degenerate one. They contain  therefore two or infinitely many. The image of the six blown-up curves $C_i$ have infinitely many.  This image moreover lies in the Clebsch cubic. There only remains the rest of the Clebsch cubic, for which there is either  one icosahedral set or infinitely many.  The whole surface is therefore the locus of sections of ${\mathcal E}^*$ in ${\bf 4}$ which have a degenerate zero set -- the second and  third cases of the trichotomy.
   
    We wish finally to determine which cubic curves contain infinitely many icosahedral sets, that is, which of the sections of ${\mathcal E}^*$ on $Z$ vanish on a positive-dimensional variety $Y\subset Z$.  Since $c_1(Z)>0$, $Y$ must intersect the anticanonical divisor $D$ nontrivially, and so we are looking for $[f]$ lying  in $S$ and also in the trisecant surface $T$. 
   
   Now the special property of  $\tilde R$ we observed is that the polar of any point on it is a tritangent plane of $R$. This means that $[f]\in \tilde R$ defines a plane section of $S\subset {\mathbf P}^3$ which meets the curve $R$ tangentially at three points. In other words the associated plane cubic $C$ passes through three and hence infinitely many icosahedral sets.  A trisecant line of $\tilde R$ meets the cubic surface $S$ generically in three points, which are the three points of intersection with $\tilde R$. Thus for these lines $T$ intersects $S$ only in $\tilde R$.
   
  There are degenerate cases --  twelve points (an $A_5$-invariant set) on $\tilde R$ where the tritangent plane to $R$ meets it in two points with multiplicity three instead of three with multiplicity two.  The trisecant line to $\tilde R$ at these points lies entirely in the Clebsch cubic  -- these are the six exceptional divisors $E_1,E_2, \dots,E_6$.  
  
  Set-theoretically, $T$ intersects $S$ in $\tilde R$ and the six lines. 
 We have therefore shown:
   
   \begin{thm} A generic point of ${\mathbf P}^3$ defines a cubic which contains precisely two icosahedral sets. If it contains one or infinitely many it lies on the Clebsch cubic surface $S$ and in the latter case it lies on the degree six  rational curve $\tilde R\subset S$  or on the six lines $E_1,\dots,E_6$. 
   \end{thm}
   
 \begin{rmk}\label{div}  In terms of  divisor classes let $H$ be the pull-back of the hyperplane divisor on the original projective plane then  a plane cubic $C$ through $[a_1],\dots,[a_6]$ 
 lifts to a curve whose divisor is  $3H-(E_1+\cdots +E_6)$. 
 The embedding for a cubic surface in ${\mathbf P}^3$ is the anticanonical embedding so
$$-K_S\sim 3H-(E_1+\cdots +E_6).$$
The divisor of $R$, the lift of a conic in ${\mathbf P}^2$, is $2H$.

Now consider the other six exceptional curves coming from the conics passing through five of the six points. These have classes
$\tilde E_1\sim 2H-(E_2+\cdots+E_6)$ etc. and 
$$-K_S\sim 3\tilde H-(\tilde E_1+\cdots +\tilde E_6)=3H-(E_1+\cdots +E_6)$$
It follows that $\tilde H \sim 5H-2(E_1+\cdots E_6)$, where $\tilde H$ is the hyperplane divisor for $\tilde {\mathbf P}^2$.
Then $\tilde R\sim 2\tilde H\sim 10H-4(E_1+\cdots E_6)$,  so the divisor class of   $T$ on $S$ is 
$$-10K_S\sim 10(3H-(E_1+\cdots +E_6))=3(10H-4(E_1+\cdots E_6))+2(E_1+\cdots E_6).$$

We may also note here that
$$R+\tilde R\sim 2H+10H-4(E_1+\cdots E_6)=4(3H-(E_1+\cdots +E_6))\sim -4K_S$$
so that there is an invariant quartic surface which vanishes on the pair of curves. This is $9\sigma_2^2-20\sigma_4$.
  \end{rmk}

   \section{The trichotomy for ${\mathbf P}(V)$}
   
Since a generic $f\in V$ defines a cubic which contains a finite number of icosahedral sets, and since  the stabilizer of an icosahedron is finite, then  ${\mathbf P}^3\subset {\mathbf P}(V)$ sweeps out an open set under the action of $SO(3,\C)$.  In principle, therefore, the trichotomy for ${\mathbf P}(V)$ entails looking for a hypersurface which intersects ${\mathbf P}^3$ in the Clebsch cubic. However, the Hilbert polynomial for  the  four-dimensional representation space  of $A_5$ is 
$$\frac{1+t^{10}}{(1-t^2)(1-t^3)(1-t^4)(1-t^{5})}$$
and for the $SO(3,\C)$ invariants in the seven-dimensional representation space 
$$\frac{1+t^{15}}{(1-t^2)(1-t^4)(1-t^6)(1-t^{10})}$$
so there is no degree $3$ invariant for $SO(3,\C)$. However, as shown in \cite{Hit}, there is an invariant   sextic hypersurface in ${\mathbf P}(V)$ which  meets ${\mathbf P}^3$ tangentially in the Clebsch cubic. The formula given in \cite{Hit} is adapted to the context of spherical harmonics and functions on the 
two-sphere, but here we shall adopt the point of view in \cite{Ig} instead, using the realization of the representation space $V$ as  $S^6$ -- homogeneous sextic polynomials in $(z_1,z_2)$. 

In inhomogeneous form the sextic is written as 
$$u_0 x^6+u_1 x^5+\dots+u_6$$
with roots $x_1,x_2,\dots,x_6$. Writing $(ij)$ for $x_i-x_j$ the discriminant is the degree $10$ invariant
$$\Delta=u_0^{10}\prod_{i<j}(ij)^2.$$
As we have seen, the vanishing of this describes the locus of cubics $f$ which are orthogonal to a degenerate isotropic subspace. In \cite{Ig}, Igusa defines similar invariants, summing over permutations:
\begin{eqnarray*}
A&=&u_0^2\sum_{\mathrm{fifteen}}(12)^2(34)^2(56)^2\\
B&=&u_0^4\sum_{\mathrm{ten}} (12)^2(23)^2(31)^2(45)^2(56)^2(64)^2\\
C&=&u_0^6\sum_{\mathrm{sixty}}(12)^2(23)^2(31)^2(45)^2(56)^2(64)^2(14)^2(25)^2(36)^2
\end{eqnarray*}
  and   introduces rational invariants $J_2,J_4,J_6,J_8,J_{10}$ where
$$J_6=\frac{A^3}{221184} + \frac{5A B}{13824} - \frac{C}{576}$$
and $J_{10}=2^{-12}\Delta$.
 \begin{thm} A generic point of ${\mathbf P}(V)$ defines a cubic which contains precisely two icosahedral sets. It contains one or infinitely many if and only if it lies on the hypersurface $J_6=0$.
  \end{thm} 
\begin{prf} The proof consists of using the explicit forms in Section 4 of \cite{Ig} to determine a degree six invariant which vanishes  on the $SO(3,\C)$ transforms of  the rational curves $\tilde R$ and $E_1,\dots,E_6$. This can be done by evaluating on the normal forms below:

\noindent (i)  A point in $\tilde R$ gives a sextic  polynomial in $z$ with three double zeros, and each such is equivalent under the action of $SO(3,\C)$ to  $z^2(z-1)^2$.   

\noindent  (ii) The cubic polynomials in  $E_i$ are of the form 
$$(c,x)((a_i,x)^2-\frac{1}{5}(a_i,a_i)(x,x))$$
with $(c,a)=0$. The line $(a,x)=0$ intersects $Q$ given by $(x,x)=0$ at two points in general and $c$ is a third point on this line. Its polar  line $(c,x)=0$ intersects $Q$ in two more points, so the sextic has two double zeros and two simple ones if $[c]$ does not lie on $Q$. There is a constraint however -- identifying $Q$ with ${\mathbf P}^1$ the four points must have cross-ratio $-1$. A normal form for this is $z^2(z^2-1)$.

We finally need to consider the points of ${\mathbf P}(V)$ whose $SO(3,\C)$ orbits do not intersect ${\mathbf P}^3$. This means cubics which do not contain any nondegenerate icosahedral set. If the cubic contains two degenerate icosahedral sets, then it is tangential twice to the null conic, so in the sextic polynomial interpretation we can take it to the form 
$$z^2(az^2+bz+c).$$
and evaluating $J_6$ gives
$$\frac{1}{1024}b^2(b^2-4ac)^2.$$
If $J_6=0$, then either $b^2-4ac=0$ which is Case (i) above or $b=0$ which is Case (ii). 
\end{prf}

\section{A special cubic curve}

We remarked above that the sextic polynomial $z^2(z-1)^2$ is a normal form for an element in $S^6$ which represents a cubic curve containing an infinite number of icosahedral sets. There is thus a single $SO(3,\C)$ orbit for a cubic which is tangential to the null conic at three distinct points. Perhaps the simplest canonical form for this is 
$$y^2=x^3+x^2+4x+4$$
where the null conic has equation
$$y^2=4(2x-3)(x+1).$$
The two meet at just three points: $(4,10),(4,-10), (-1,0)$.

Note that this cubic $C$  is defined over the integers and has rational points $$\infty,(4,-10),(0,-2),(-1,0),(0,2),(4,10)$$ forming a cyclic group of order six.

Now $C$ contains a one-parameter family of icosahedral sets, three of which are degenerate. Moreover this  family is rational. To see this, we revert to the isotropic subspace description.  From Proposition \ref{pencil} we may assume that there are two isotropic subspaces 
$U_1,U_2$  of $W$ for which $U_1\cap U_2=0$.   A third subspace $U_3$ intersects each of these trivially and hence is the graph of a  linear transformation $S:U_1\rightarrow U_2$. Thus for $w_1,w_2\in U_1$ we have
$$0=\omega_x(w_1+Sw_1,w_2+Sw_2)=\omega_x(w_1,Sw_2)+\omega_x(Sw_1,w_2)$$
since $\omega_u(w_1,w_2)=0=\omega_u(Sw_1,Sw_2)$. This means the graph of $tS$ for any  $t\in \C$ is also  isotropic, so we have a  family of isotropic subspaces with rational parameter $t$.

Each icosahedral set consists of six points and so we have a sixfold covering $p:C\rightarrow {\mathbf P}^1$ with three branch points $0,1,\infty$ in ${\mathbf P}^1$. At these points the icosahedron degenerates -- one axis remains (the vector $b$ in the type 1 degenerate subspace) and the other five  coalesce to be a null vector (the vector $a$). The inverse image of each branch point consist of two points -- one, $[a]$, a ramification point of order $5$ and one simple point $[b]$. These six points are precisely the six rational points in the canonical form above. 

The projection map $p$ is given in its simplest form by the meromorphic function
$$p^2=\frac{(1+x)^5(y-x-2)}{(y+3x-2)^5}.$$
An explicit form for $p$ is a little less simple -- it involves elliptic functions for  two planar embeddings (see Remark 5 in Section 12), which differ by  the line bundle $L$ of order $2$ , for here we are considering a curve  where the vector bundle is $L\oplus L$ with $L^2$ trivial.

\begin{rmk} The formulae above have come from \cite{Vid} in the context of algebraic solutions to the hypergeometric equation. The monodromy of the covering $p$ lies in the icosahedral group and mapping into $PSL(2,\C)$ this gives a second order differential equation with three singular points. Schwarz's celebrated paper \cite{S} classifying these stops short of explicit expressions and the elliptic curve above appears in an attempt to be more concrete about these. 
\end{rmk}

  \section{Higher degree curves}
  The picture for higher degree representations and curves is not clear -- in particular there are stable bundles and not  just semi-stable as in the case of the elliptic curve. One might ask for a link between semistability and isotropic subspaces, however. The argument in Section \ref{iso} generalizes to show that a $d$-dimensional  isotropic subspace orthogonal to a harmonic polynomial $f$ of degree $d$ defines a line subbundle $L$ with $\deg L=\deg E/2$, though it is not clear if the reverse is true. What is  true is that isotropic subspaces are rare. We can find some, though, as we show below.

Let $V$ now be the $(2d+1)$-dimensional irreducible representation space of $SO(3,\C)$ and consider a $d$-dimensional subspace isotropic for all skew forms $(x\cdot u,v)$ which is invariant by  $SO(2,\C)\subset  SO(3,\C)$. Then it is a sum of $d$ weight spaces $L_m$ where $-d\le m\le d$. Let ${\mathbf h},{\mathbf n}^+,{\mathbf n}^-$ be a basis of ${\lie {so}}(3)$ where ${\mathbf h}$ generates $SO(2,\C)$ and $[{\mathbf h},{\mathbf n}^+]={\mathbf n}^+,  [{\mathbf h},{\mathbf n}^-]=-{\mathbf n}^-$.
\begin{prp} Let $U\subset V$ be a $d$-dimensional  $SO(2,\C)$-invariant isotropic subspace of $V$. Then $U$ is a sum of weight spaces $\pm  \{1,2,3,\dots, d\}$ or $\pm \{0,2,\dots,d\}$
\end{prp}
\begin{prf}
Because ${\mathbf h}$ acts the  scalar $m$ on $L_m$, if a weight $m>0$ occurs, then $-m$ does not  because on $L_m\oplus L_{-m}$, $({\mathbf h}\cdot u,v)=muv$ which is non-zero. 

Now ${\mathbf n}^+:L_m\mapsto L_{m+1}$, so on $L_m\oplus L_{-m-1}$,  $({\mathbf n}^+\cdot u,v)=uv$ which is non-zero. Hence if $m$ occurs, $-m-1$ does not. Similarly with ${\mathbf n}^-$, if $m$ occurs $-m+1$ does not. 

Suppose $0$ is not one of the $d$ weights and let $P$ be the positive weights. Then the negative weights $N$ must be the negatives of the complement of $P$ in $\{1,\dots,d\}$ by the first criterion. If $m\in P$ then by the second criterion, the adjacent numbers $m-1,m+1$ must be in $P$. But this means that $P$ consists of $1,2,3,\dots, d$. If $P$ is empty, then we get $-1,-2,-3,\dots, -d$

If $0$ is one of the weights, then by the second criterion $\pm 1$ is not a weight. Then $P$ is a subset of $2,\dots,d$. Again the $d-1$ non-zero weights must be made up of $P$ and the negative of its complement. As before adjacent positive numbers  must be in $P$, so $P$ is $2,3,\dots, d$ so the weights are 
$0,2,\dots,d$ or their negatives. 
\end{prf}
Note that these are precisely the {\it degenerate} isotropic subspaces which we have encountered for $d=3$. So if the curve $f(x)=0$ has a high order contact with the conic $Q$, we may deduce that the bundle $E$ is strictly semistable.

\section{Further remarks}

\noindent 1. We have shown here that a generic cubic which passes through one icosahedral set passes through another one. According to Melliez \cite{Mel} (who also introduced the use of the bundle ${\mathcal E}^*$), this is a result of Reye and can be found in Baker \cite{Bak}, page 145, Exercise 26 or Coble \cite{Cob} page 236. The reader is invited to make the translation.

\noindent 2. A discussion of the real case  can be found in \cite{Hit}.  There we consider the picture in the sphere double covering $\RP^2$. The inverse image of the cubic curve is then the so-called {\it nodal set} of a spherical harmonic and the issue is whether the nodal set contains the vertices of a regular icosahedron.The rational curves $R$ and $\tilde R$ are not real and so do not appear in the story. 

\noindent 3. 
The curve $\tilde R$ is simply the null conic in $\tilde{\mathbf P}^2$ but is much more complicated when viewed in ${\mathbf P}^2$. We saw from Remark \ref{div} that its divisor class is $10H-4(E_1+\cdots E_6)$ so It is of degree $10$.  In fact it is the  invariant singular  curve defined by: 
$ 2({{x_1}}^{10} +  {{x_3}}^{10}+{{x_2}}^{10}) 
  +35({{x_1}}^8{{x_2}}^2 +{{x_1}}^2{{x_2}}^8+  
  {{x_2}}^2{{x_3}}^8 + 
  {{x_2}}^8{{x_3}}^2+ 
  {{x_3}}^8{{x_1}}^2+{{x_3}}^2{{x_1}}^8)+ 
  25{\sqrt{5}}({{x_1}}^2
   {{x_2}}^8+ 
  {{x_2}}^2
   {{x_3}}^8 +{{x_3}}^2
  {{x_1}}^8
   - {{x_1}}^8
   {{x_2}}^2- 
  {{x_2}}^8
   {{x_3}}^2- 
  {{x_3}}^8{{x_1}}^2
   )
  -30({{x_1}}^6{{x_2}}^4 + 
  {{x_1}}^4{{x_2}}^6+ 
   {{x_2}}^6{{x_3}}^4+ {{x_2}}^4{{x_3}}^6+{{x_3}}^6{{x_1}}^4
   + 
  {{x_3}}^4{{x_1}}^6 
  )+
  50{\sqrt{5}}({{x_1}}^6
   {{x_2}}^4 +   
  {{x_2}}^6
   {{x_3}}^4+{{x_3}}^6 
 {{x_1}}^4
   -{{x_1}}^4
   {{x_2}}^6 
     - 
  {{x_2}}^4
   {{x_3}}^6-{{x_3}}^4
  {{x_1}}^6
   ) 
  -  560x_1^2x_2^2x_3^2(x_1^4+$ $+x_2^4+x_3^4)
 +
  1060 x_1^2x_2^2x_3^2(x_1^2x_2^2+x_2^2x_3^2+x_3^2x_1^2)=0.$

\noindent 4. The closure of the curves in the Mukai-Umemura threefold defined by icosahedra with a  fixed axis are rational curves $C$ with $c_1(Z)[C]=2$, i.e degree $2$. These are the minimal rational curves appearing in the classification of Fano varieties. The one-parameter family of icosahedra parametrized by points of $\tilde R$ define degree $3$ curves, for the plane cubic $C$ is tangential to the null conic at three points. Whereas degree $3$ rational curves were linked to icosahedral  solutions to the hypergeometric equation, degree $4$ curves give rise to algebraic solutions of a particular case of Painlev\'e's sixth equation. This is the setting of \cite{Hit1},\cite{Hit2} where the dihedral and octahedral cases are examined.  A classification of all such solutions, akin to Schwarz's list, has been given by Boalch \cite{Boalch}.

To find degree $4$ rational curves we consider the map $\beta:\hat Z\rightarrow {\mathbf P}^3$ of Section \ref{fromMuk}. Since $\tilde R$ parametrizes cubics containing an infinite number of icosahedral sets, the map $\beta$ collapses a subvariety in   $\hat Z$. If we want to go the other way, from ${\mathbf P}^3$ to $Z$, we first blow up $\tilde R$ and then do some blowing up and down. Now take a {\it secant} to the curve $\tilde R\subset {\mathbf P}^3$ -- a line joining two points. A generic line intersects   the trisecant surface $T$ in  $10$ points. However, the $3:3$ correspondence between $R$ and $\tilde R$ implies that through a generic point of  $\tilde R$ there pass three trisecants. Blowing up  $\tilde R$ means that the proper transform of  the secant is a curve which meets the blown-up trisecant surface in $10-3-3=4$ points, and this becomes a rational curve of degree $4$ in the Mukai-Umemura threefold. The corresponding solution to Painlev\'e VI is known and is due to Dubrovin and Mazzocco \cite{DM}. An open set in the Mukai-Umemura threefold can then be identified with the twistor space of a complex self-dual four-manifold and the secants above, after transforming by $SO(3,\C)$, give the four-parameter family of twistor lines.

\noindent 5. One might ask, without reference to the vector bundle $E$,  how a  cubic curve defined by a harmonic polynomial gives rise to a  line bundle $L$. It arises from the geometry as follows.

As $a$ varies in the cubic $C$ consider the polynomial $f_a$, invariant under the action of $a\in {\lie{so}}(3)$. Now $\omega_a(f_a,v)=(a\cdot f_a,v)=0$ so that $f_a$ generates a line bundle in the kernel of $\omega: W\rightarrow W^*(1)$. We identify this bundle as follows. Consider the map
$\gamma:{\mathbf P}^2\rightarrow {\mathbf P}(V)$ defined by $a\mapsto f_a$. Then $f_a$ spans the pull-back $\gamma^*({\mathcal O}(-1))$ and since $f_a$  is homogeneous of degree $3$ in $a$, $\gamma^*({\mathcal O}(-1))\cong{\mathcal O}_{{\mathbf P}^2}(-3)$. Restrict this to  $C$ to get the required line bundle.

We then have an inclusion 
${\mathcal O}(-3)\subset E^*(-1)$ over $C$. Projecting to  $E^*(-1)/L(-1)$ gives a section of $L^*(2)$. This vanishes when $f_a\in U$ which  is at the six points $[a_1], \dots, [a_6]$ defined by the icosahedral set. Thus on $C$ we have the relation of divisor classes 
$$L\sim {\mathcal O}(2)- \sum_1^6[a_i].$$
Now we saw that the cubics through the icosahedral set are given by plane sections of the Clebsch cubic surface. This means that the same cubic is embedded in two different planes -- a plane in ${\mathbf P}^3$ and the original plane with the conic $Q$. Since the embedding in ${\mathbf P}^3$ is 
$${\mathcal O}(3)- \sum_1^6[a_i]$$
we see that the degree zero line bundle $L$ is the {\it difference} of the two hyperplane divisor classes.

\vskip 1cm
 Mathematical Institute, 24-29 St Giles, Oxford OX1 3LB, UK
 
 hitchin@maths.ox.ac.uk

 \end{document}